\documentclass[11pt]{article}
\usepackage{amssymb}
\usepackage{amsmath}
\usepackage{amsfonts}
\usepackage{amsthm}

\numberwithin{equation}{section}

\newtheorem{theorem}{Theorem}[section]

\newtheorem{proposition}{Proposition}[section]

\def\R{\mathbb{R}}

\def\N{\mathbb{N}}

\def\T{\mathbb{T}}

\begin{document}

\title{A note on the inviscid limit of the Benjamin-Ono-Burgers equation in the energy space}
\author{\normalsize \bf  Luc Molinet  \\
{\footnotesize \it L.M.P.T., Universit\'e Fran\c cois Rabelais Tours, F\'ed\'eration Denis Poisson-CNRS,} \\
{\footnotesize \it Parc Grandmont, 37200 Tours,  France} }
\maketitle


\vspace*{1mm} \vspace*{2mm}\\

\begin{abstract}
In this paper we study the inviscid limit of the Benjamin-Ono-Burgers equation in the energy space $ H^{1/2} (\R) $ or $ H^{1/2}(\T) $. We prove the strong convergence in the energy space of the solution to this equation toward the solution of the Benjamin-Ono equation as the dissipation coefficient converges to $ 0 $.

\end{abstract}

\maketitle

\section{Introduction}

The initial value problem (IVP) associated to the Benjamin-Ono-Burgers equation
\begin{equation}\label{BOB}
\left\{\begin{array}[pos]{ll}
           \partial_tu+\mathcal{H}\partial^2_xu-\varepsilon u_{xx}=u\partial_xu\\
           u(x,0)=u_0(x),
       \end{array} \right.
\end{equation}
where $x\in \R $ or $ \T $, $t \in \mathbb R$, $u$ is a real-valued function, $ \varepsilon $ is a positive real number and
$\mathcal{H}$ is the Hilbert transform given for function on $ \R $ by 
\begin{equation} \label{hilbert}
\mathcal{H}f(x)=\text{p.v.}\,\frac1\pi\int_{\mathbb  
R}\frac{f(y)}{x-y}dy.
\end{equation}
This equation is a dissipative perturbation of the celebrated Benjamin-Ono equation and is proven to be physically relevant in \cite{ER}. Recall that the Benjamin-Ono equation   was derived in \cite{B} to modelize   the unidirectional evolution of weakly nonlinear dispersive internal
long waves at the interface of a two-layer system, one being
infinitely deep.  

It is easy to check that the BOB equation is globally well-posed in $ H^s(\R) $ and $ H^s(\T) $  for $ s>-1/2 $. One can use for instance purely parabolic approach as for the dissipative Burgers equation (cf. \cite{Be}). On the other hand, the I.V.P. for the BO equation is more delicate to solve in Sobolev spaces with low indexes. Indeed, it was proven in \cite{MST} that the this I.V.P. cannot be solved by a Picard iteration scheme in any Sobolev space $ H^s(\R) $. However,  in \cite{Ta} Tao noticed that one can solve this I.V.P. in $ H^1(\R) $ by introducing a suitable gauge transform. This approach has been pushed forward in \cite{BP}, \cite{IK} and \cite{MP}. In these two last papers, the I.V.P. is proved to be globally well-posed in $ L^2(\R) $ (see  \cite{L} or \cite{MP}  for the global well-posedness in $ L^2(\T)$). Unfortunately this gauge transform does not behave well with respect to perturbations of the equation and in particular with respect to the BOB equation. It is thus not evident to prove the uniform in $ \varepsilon $ wellposedness of the BOB equation in low regularity spaces and , as a consequence, to study its inviscid limit behavior. However, in \cite{guo}, Guo {\it  and all}  used the variant of the Bourgain's spaces introduced in \cite{KT} to prove that the BOB equation is uniformly (in $\varepsilon$) well-posed in $ H^s(\R) $, $s\ge 1$, and to deduce the convergence of the solutions to this equation towards the one of the BO equation in $ H^s(\R) $, $s\ge 1 $.

 Our goal in this paper is to prove that this convergence result also holds in the energy space $ H^{1/2}(\R) $ by a very much simpler approach. This approach combines the conservation laws  and the unconditional uniqueness  in $H^{1/2} $  (cf. \cite{BP}, \cite{MP}) of the Benjamin-Ono equation.   Note that this approach works also in  $H^{1} $ (and more generally in all $ H^{n/2}$ for $ n\ge 1$) where the unconditional uniqueness is a simple consequence of the $L^2 $ Lipschitz bound established in 
\cite{Ta}. Therefore, our approach gives also a great simplification of the inviscid limit result in $ H^1(\R) $ with only \cite{Ta} in hand. Finally, it is worth noticing that our method also works in the periodic setting.
\section{Main result and proof}
Our main result reads 
\begin{theorem}
Let $ K: =\R $ or $ \T $, $ u_0\in H^{1/2}(K) $ and $ \{\varepsilon_n\}_{n\in\N} $ be a decreasing sequence of real numbers converging to $ 0 $. 
For any $ T>0 $, the sequence  
$ \{u_n\}_{n\in\N} \in C(\R_+;H^{1/2}(K))$   of solution to \eqref{BOB}$_{\varepsilon_n}$ emanating from $ u_0 $ satisfies
 \begin{equation}\label{to1}
u_{n} \to u \mbox{ in }C([0,T];H^{1/2}(K))
\end{equation}
where $ u\in C(\R;H^{1/2}(K)) $ is the unique solution to the BO equation emanating from $u_0$.
 \end{theorem}
 \proof We first give the complete proof in the real line case and discuss the adaptation in the periodic case at the end of the paper. We will divide the proof in three steps. \vspace*{2mm}\\
{\bf Step 1.} {\it Uniform bounds.} \\
 First we establish  uniform in $ \varepsilon $ a priori estimates on solutions to \eqref{BOB}$_\varepsilon $.
Taking the $ L^2 $-scalar product of the equation with $ u$ it is straightforward to check  that smooth solutions to \eqref{BOB}$_{\varepsilon_n}$ satisfy 
$$
\frac{d}{dt} \|u_n (t)\|_{L^2_x}^2+\varepsilon \|D_x u_{n} (t) \|_{L^2_x}^2 =0 
$$
By the continuity  of the flow-map of  \eqref{BOB}$_{\varepsilon_n}$ in $ H^{1/2}(\R) $ it follows that for any $ u_0\in H^{1/2} $, 
 $ u_n\in C_b(\R_+;L^2(\R))$ and $ u_{n,x}\in  L^2 (\R_+;H^1(\R)) $ with 
\begin{equation}\label{toto}
\|u_n (t)\|_{L^2_x}^2+\varepsilon \int_0^\infty \|\partial_x u_n(s)\|_{L^2_x}^2  \, ds  \lesssim \|u_0\|_{L^2}^2
\end{equation}
Similarly, taking the $ L^2 $-scalar product of the equation with $ D_x u_n +u^2_n/2 $ and setting 
$$
E(v):=\frac{1}{2} \int_{\R} |D_x^{1/2} v |\, dx + \frac{1}{6} \int_{\R} v^3 \, dx 
$$
 we obtain 
\begin{eqnarray*}
 \frac{d}{dt} E(u_n(t)) +\varepsilon \|D_x^{3/2}  u_{n} (t) \|_{L^2_x}^2  &= &  \frac{\varepsilon}{2} \int_{\R} u^2_n u_{n,xx} \\
 & \lesssim & \varepsilon \|u_n(t)\|_{L^\infty_x} \|u_{n,x}\|_{L^2_x}^2 \\
 & \lesssim & \varepsilon \|u_n(t)\|_{L^\infty_x} \|D_x^{1/2} u_n\|_{L^2_x}\|D^{3/2}_x u_n\|_{L^2_x}
\end{eqnarray*}
Therefore by Young inequality and then classical interpolation inequalities, we infer that 
\begin{equation}\label{tata}
 \frac{d}{dt} E(u_n(t)) +\frac{\varepsilon}{2} \|D_x^{3/2}  u_{n} (t) \|_{L^2_x}^2 \lesssim \varepsilon \|u_n(t)\|_{L^\infty_x}^2 \|D_x^{1/2} u_n\|_{L^2_x}^2
\lesssim \varepsilon \|u_n(t)\|_{L^2_x}^2 \| u_{n,x}\|_{L^2_x}^2
\end{equation}
Gathering this last estimate with \eqref{toto} we obtain 
$$
 E(u_n(t))+\frac{\varepsilon}{2} \int_0^t \|D_x^{3/2}  u_{n} (s) \|_{L^2_x}^2\, ds \lesssim \|u_0\|_{L^2_x}^4+E(u_0), \quad \forall t\ge 0 \; .
 $$
Using classical Sobolev inequalities and again \eqref{toto} this eventually leads to 
$$
 \|D_x^{1/2} u_n (t)\|_{L^2_x}^2 +\varepsilon \int_0^t \|D_x^{3/2}  u_{n} (s) \|_{L^2_x}^2\, ds \lesssim \|u_0\|_{H^{1/2}_x}^4 +\|u_0\|_{L^2_x}^2 \Bigl(  \|D_x^{1/2} u_0\|_{L^2_x}+ \|D_x^{1/2} u_n (t)\|_{L^2_x}\Bigr) , \quad \forall t\ge 0 \; .
 $$
which ensures that 
\begin{equation}\label{titi}
 \|D_x^{1/2} u_n (t)\|_{L^2_x}^2 +\varepsilon \int_0^t \|D_x^{3/2}  u_{n} (s) \|_{L^2_x}^2\, ds \lesssim 1+\|u_0\|_{H^{1/2}_x}^4, \quad \forall t\ge 0 \; .
\end{equation}
{\bf Step 2.} {\it Convergence in the weak topology.}
\begin{proposition}\label{pro1}
Let $ u_0\in H^{1/2}(\R) $ and $ \{\varepsilon_n\}_{n\in\N} $ be a decreasing sequence of real numbers converging to $ 0 $. 
The sequence  
$ \{u_n\}_{n\in\N} $   of solution to \eqref{BOB}$_{\varepsilon_n}$ emanating from $ u_0 $ satisfies
 \begin{equation}\label{to1}
u_{n} \rightharpoonup u \mbox{ weak star in } L^\infty(\R;H^{1/2}(\R))
   \end{equation}
and 
   \begin{equation} \label{to2}
   (u_{n}, \phi)_{H^{1/2}} \to (u,\phi)_{H^{1/2}} \mbox{ in } C([-T,T]) , \forall \phi\in H^{1/2}(\R),\;   \forall  T>0.
   \end{equation}
where $ u\in C(\R;H^{1/2}(\R)) $ is the unique solution to the BO equation emanating from $u_0$.
\end{proposition}
\proof 
According to \eqref{toto} and \eqref{titi}, the sequence $\{u_n\} $ is   bounded in  $L^\infty(\R;H^{1/2}(\R))$. Moreover, in view of the equation \eqref{BOB},
 the sequence $ \{\partial_t u_{n}\}$ is bounded in $ L^\infty(\R;H^{-2}(\R))$. By Aubin-Lions compactness theorem, we infer that for any $ T>0 $ and $ R>0 $,  $\{u_n\} $ is relatively compact in $ L^2(]-T,T[\times]-R,R[) $.Therefore, using a diagonal extraction argument, we obtain the existence  of an increasing sequence $\{n_k\}\subset \N $ and $ u\in L^\infty(\R;H^{1/2}(\R)) $ 
  such that 
  \begin{align}
   u_{n_k} \rightharpoonup u \mbox{ weak star in } L^\infty(\R;H^{1/2}(\R)) \label{conv1} \\
   u_{n_k} \to u \mbox{ in }L^2_{loc}(\R^2)  \label{conv2}\\
   u_{n_k} \to u \mbox{ a.e. in }  \R^2  \label{conv3} \\
      u_{n_k}^2\rightharpoonup  u^2 \mbox{ weak star in } L^\infty(\R;L^2(\R)) \label{conv4} 
   \end{align}
 These convergences results enable us to pass to the limit on the equation and to obtain that the limit function $ u$ satisfies 
  the BO equation in the distributional sense. Now, the crucial argument is that, according to \cite{MP}, BO is unconditionally well-posed in $ H^{1/2}(\R) $(and even in $H^{1/4}(\R) $) in the sense that the solution constructed in \cite{MP} is the only function belonging to 
   $ L^\infty(-T,T;H^{1/2}(\R)) $  that satisfies (BO) in the distributional sense and is equal\footnote{Note that according to the equation,  the time derivative of such a function  must belong to $ L^\infty(-T,T; H^{-2}) $ and thus such function has to belong to
   $ C(-T,T;H^{-2})$} to $ \varphi$ at $ t=0$.  By the uniqueness of the possible limit, this ensures that, actually,
    $\{u_{n}\} $ converges to $ u $ in the sense \eqref{conv1}-\eqref{conv4}.
  
  Now, using the equation and the bound on $ \{u_n\} $, it is clear that for any $ \phi\in C^\infty_0(\R) $ and any $ T>0 $, the sequence 
   $\{t\mapsto (u_n, \phi)_{H^{1/2}}\} $ is uniformly equi-continuous on $ [-T,T] $. By Ascoli's theorem it follows  that 
   $\{ (u_{n}, \phi)_{H^{1/2}}\} \to (u,\phi)$ in $C([-T,T]$. Since  $\{u_n\} $ is bounded in $L^\infty(\R;H^{1/2}(\R))$,
    this yields \eqref{to2}.  \vspace*{2mm}  \\
  {\bf Step 3.} {\it Making use of the conservation laws of the BO equation}  \\
We start by  proving a strong convergence result in $ L^2$. 
A first idea would be to derive  a Lipschitz bound in $ L^2$. Note that, following \cite{Ta}, such a  $ L^2$-Lipschitz bound would be  available at the $ H^1 $-regularity.  However, we do not know how to get it at the $H^{1/2}$-regularity. Instead, we will rely on the non increasing property  of the $ L^2 $-norm of solutions to \eqref{BOB}$_{\varepsilon_n}$ established in \eqref{toto}. According to the $ L^2 $-conservation law of the BO equation, this proves that for any $ n\in \N $ and $ t> 0 $,
$$
 \| u_n(t)\|_{L^2} \le \|u_0\|_{L^2}=\|u(t)\|_{L^2} \; .
$$
This together with \eqref{to2} ensures that for any $ T> 0 $, 
\begin{equation}\label{zs}
u_n \to u \mbox{ in } C([0,T]; L^2(\R))  .
\end{equation}
   Combining this last convergence resultl with Proposition \ref{pro1} we infer that for any $ 0<s<1/2$, and any $ T> 0 $, 
    \begin{equation}\label{to6}
u_{n} \to u \mbox{ in } L^\infty(0,T;H^s(\R)) 
   \end{equation}
   Now, from  interpolation inequalities, \eqref{toto} and \eqref{titi} we infer that for any fixed $ t>0 $, 
   \begin{eqnarray*}
   \varepsilon \int_0^t \|u_n(s) \|_{L^2_x}^2  \|\partial_x u_n(s) \|_{L^2_x}^2\, ds &  \lesssim  & 
    \varepsilon \int_0^t \|u_n(s) \|_{L^2_x}^{8/3}   \|D_x^{3/2}  u_n(s) \|_{L^2_x}^{4/3}\, ds \\
   &   \lesssim  &  \|u_0\|_{L^2_x}^{8/3}\varepsilon^{1/3} t^{1/3} \Bigl[\varepsilon \int_0^t  \|D_x^{3/2}  u_n(s) \|_{L^2_x}^2 \, ds \Bigr]^{2/3}\\
   & \lesssim & O(\varepsilon^{1/3}) \; .
      \end{eqnarray*}
    We thus deduce from \eqref{tata} that for any $ t\ge 0 $, 
    $$
    E(u_n(t)) \le E(u_0)+O(\varepsilon^{1/3})\longrightarrow E(u_0)=E(u(t))\;  .
 $$
  Moreover, according to \eqref{to6}, , 
  $$
  \sup_{t\in [-T,T]}\Bigl|\int_{\R} u_n^3(t,x) \, dx -\int_{\R} u^3(t,x) \, dx \Bigr| \to 0  \mbox{ as } n\to \infty 
  $$
  Gathering the two above convergence results, we deduce that for all $ t\in \R $, 
  $$
  \limsup_{n\to \infty} \int_{\R} |D^{1/2}_x u_n(t)  |^2\, dx \le   \int_{\R} |D^{1/2}_x u(t)  |^2\, dx 
  $$
  which combined  with  \eqref{to2} and \eqref{zs}  ensures that for any $ T>0 $, 
 $$
 u_n \to u \mbox{ in } C([0,T]; H^{1/2}(\R)) \: .
 $$
This concludes the proof of the theorem when $ K=\R $.\vspace{2mm} \\
Finally it is easy to check that the case $ K=\T $ works exactly as well since according to \cite{MP}  the Benjamin-Ono equation is unconditionnaly well-posed in $ H^{1/2}(\T) $. Actually, this case is even simpler since the weak convergence in $H^{1/2}(\T) $ directly implies the strong convergence in $ H^s(\T ) $ for $ s<1/2$. 
\bibliographystyle{amsplain}

\providecommand{\bysame}{\leavevmode\hbox to3em{\hrulefill}\thinspace}
\providecommand{\MR}{\relax\ifhmode\unskip\space\fi MR }
\providecommand{\MRhref}[2]{%
  \href{http://www.ams.org/mathscinet-getitem?mr=#1}{#2}
}
\providecommand{\href}[2]{#2}
\begin{thebibliography}{}

\end{thebibliography}


\begin{thebibliography}{99}


\bibitem {Be} 
\newblock  D.Bekiranov,
\newblock \emph{The initial-value problem for the generalized Burgers' equation}, 
\newblock  Diff. Int. Eq., {\bf 9}  (1996), pp. 1253--1265.

\bibitem{B} 
\newblock  T.B. Benjamin,
\newblock \emph{Internal waves of permanent form in fluids of great depth,}
\newblock  J. Fluid Mech. {\bf 29} (1967), 559--592.

\bibitem{BP} 
\newblock N. Burq and F. Planchon,
\newblock \emph{On well-posedness for the Benjamin-Ono equation,}
\newblock Math. Ann., \textbf{340} (2008), 497--542.

\bibitem{BP2} N. Burq and F. Planchon,
\newblock \emph{The Benjamin-Ono equation in energy space,}
\newblock Phase
space analysis of partial differential equations, 55--62,
Progr. Nonlinear Differential Equations Appl.,  \textbf{69}, Birkh{\"a}user Boston,
Boston, MA, 2006.

\bibitem{ER}
\newblock P.M. Edwin, B. Roberts,
\newblock \emph{The Benjamin-Ono-Burgers equation: an application in solar physics,}
\newblock Wave Motion 8 (2) (1986) 151Ð158. 


\bibitem{guo}
\newblock Z. Guo, L. Peng, B. Wang and  Y. Wang
\newblock \emph{Uniform well-posedness and inviscid limit for the Benjamin-Ono-Burgers equation}  ,
\newblock Advances in Mathematics \textbf{228} (2011),  647-677.

\bibitem{IK}
\newblock A. Ionescu and C. E. Kenig,
\newblock \emph{Global well-posedness of the Benjamin-Ono equation in low-regularity spaces,}
\newblock  J. Amer. Math. Soc., \textbf{20} (2007), 753--798.


\bibitem{KT} H. Koch and D.Tataru,
 \newblock\emph{A priori bounds for the 1D cubic NLS in negative Sobolev spaces}
 \newblock  Int. Math. Res. Not. (2007), no. 16.
 
\bibitem{L} L. Molinet, 
\newblock \emph{Global well-posedness in $L^2$ for the periodic Benjamin-Ono equation,} 
\newblock Amer. J. Math., \textbf{130} (2008), 635--683. 
\bibitem{MP} L. Molinet and D. Pilod, 
\newblock \emph{The Cauchy problem for the Benjamin-Ono equation in $L^2$ revisited,} 
\newblock to appear in Anal. PDE.
\bibitem{MST} L. Molinet, J-C. Saut and N. Tzvetkov,
\newblock \emph{Ill-posedness issues for the Benjamin-Ono and related equations,}
\newblock SIAM J. Math. Anal.,  \textbf{33} (2001), no. 4, 982--988.

\bibitem{Ta}
\newblock T. Tao,
\newblock \emph{Global well-posedness of the Benjamin-Ono equation in $H^1(\mathbb R)$,}
\newblock J. Hyp. Diff. Eq. \textbf{1} (2004), 27--49.

\end{thebibliography}

\end{document}